\documentclass[12pt]{article}

\usepackage{latexsym,amsfonts,amsthm,amsmath,amscd,amssymb}
\usepackage[dvips]{graphicx}

\newtheorem{lemma}{Lemma}[section]
\newtheorem{thm}[lemma]{Theorem}
\newtheorem{rem}[lemma]{Remark}
\newtheorem{prop}[lemma]{Proposition}

\newtheorem{defn}[lemma]{Definition}

\def\sphere{{\mathbb S}}

\def\C{{\mathbb C}}
\def\R{{\mathbb R}}
\def\H{{\mathbb H}}
\def\Z{{\mathbb Z}}
\def\N{{\mathbb N}}

\def\D{{\cal D}}
\def\E{{\cal E}}
\def\H{{\cal H}}
\def\G{{\cal G}}
\def\T{{\mathbb T}}
\def\qsfera{{\mathbb S}^2_{q,c}}
\def\O{{\cal O}}

\def\g{\mathfrak{g}}

\def\h{{\mathfrak h}}

\def\mG{\mathcal{G}}

\def\mP{\mathcal{P}}

\begin{document}

\title{The quantization of the symplectic groupoid of the standard Podle\`s sphere}

\author{F. Bonechi\footnote{\small INFN Sezione di Firenze, email: bonechi@fi.infn.it} , 
N. Ciccoli\footnote{\small Dipartimento di Matematica, Univ. di Perugia, email: ciccoli@dipmat.unipg.it} ,
N. Staffolani\footnote{\small INFN Sezione di Firenze, email: nicola.staffolani@gmail.com} ~and 
M. Tarlini\footnote{\small INFN Sezione di Firenze,  email: marco.tarlini@fi.infn.it}
}

\date{September 19, 2012}

\maketitle

\begin{abstract}
We give an explicit form of the symplectic groupoid $\G(\sphere^2,\pi)$ that 
integrates the semiclassical standard Podle\`s sphere ($\sphere^2,\pi$). 
We show that Sheu's groupoid $\G_S$, whose convolution $C^*$-algebra quantizes 
the sphere, appears as the groupoid of the Bohr-Sommerfeld leaves of a 
(singular) real polarization of $\G(\sphere^2,\pi)$. By using a complex polarization 
we recover the convolution algebra on the space of polarized sections. We stress the 
role of the modular class in the definition of the scalar product in order to get the 
correct quantum space.
\end{abstract}

\thispagestyle{empty}

\eject

\section{Introduction}
Quantization, very generally, can be seen as a process which aims at replacing a set 
of commutative classical observables with non commuting quantum ones. There are two 
distinct ways to achieve this goal: either changing the nature of observables, from
functions to linear operators, or changing the product from pointwise product to a 
non commutative $\ast$-product. Roughly speaking these two approaches correspond 
respectively to \emph{geometric quantization} and \emph{deformation quantization}.
Integration  by symplectic groupoids was introduced, from the beginning, as a tool to 
quantize Poisson manifolds, {\it i.e.} when the Poisson bracket is degenerate. Generally
speaking the quantization of the graph of the groupoid multiplication should equip the 
Hilbert space of states with a $*$-algebra structure, which can be seen as the
quantization of the Poisson manifold of units. The geometric prequantization procedure 
was fully developed in \cite{WeXu}, where it is shown that it can be interpreted as a 
groupoid extension. More difficult problems, in this context, come into play when one has 
to deal with polarizations, which has been largely ignored since then. The full picture of 
quantization through symplectic groupoids was clarified only in few cases, the
main one being the non commutative irrational torus in \cite{We934}.

Recently in \cite{Hawkins} it was introduced a notion of multiplicative polarization taking 
into account compatibility conditions between  the symplectic and groupoid structure. In the
best possible case the set of Lagrangian leaves of the polarization inherits a structure of 
topological groupoid itself. It is then possible to consider the convolution $C^*$-algebra of
this groupoid as being a quantization of the original Poisson manifold. The difficulties here 
reflects the usual rigidity of geometric quantization: how to deal with compact leaves and 
consequent Bohr-Sommerfeld rules, need of complex polarizations, etc....

Quantum groups and related homogeneous spaces are examples which have been extensively studied 
at the quantum level and a lot is known about the $C^*$-algebras of these quantum spaces. In
particular in a series of papers A. Sheu \cite{Sheu1,Sheu3} has observed that $C^*$-algebras of 
quantum groups are related to specific groupoid $C^*$-algebras. In his approach such groupoids
appear from a detailed analysis of the combinatorics of the space of irreducible $*$--representations 
of the $C^*$--algebra (or from the symplectic foliation of the underlying Poisson manifold - the
two being strictly linked). We want to investigate how these groupoids can be related to the geometric 
quantization of the symplectic groupoid integrating the semiclassical Poisson manifold. In this paper 
we will restrict this analysis to the simplest example of quantum homogeneous space, the so called 
{\it standard Podle\`s sphere}. The $C^*$-algebra of this quantum homogeneous space is simply the 
unitization of compact operators, as observed already in \cite{Podles}; according to \cite{Sheu3}, 
it can be described as the groupoid $C^*$-algebra of a particular subgroupoid $\G_S$ of the so 
called Cuntz groupoid $\O_1$. This $C^*$--algebra can be seen as a quantization of a $SU(2)$--covariant 
Poisson structure on the sphere having only two symplectic leaves: a degenerate point with
quadratic singularity and a symplectic plane (see \cite{Sheu2}). Despite its apparent simplicity, this example 
is already quite involved from the point of view of quantization: for instance the Poisson cohomology 
class of the Poisson tensor is not trivial so that the prequantization cocycle is not trivial too. 
Moreover, the modular class is not zero, so that a modular operator is expected to come out in the 
quantization.

The first step of the program is to construct an explicit description of the symplectic groupoid integrating 
the semiclassical sphere and to discuss the modular function as the function that integrates the modular 
vector field. We will use a description of the integration of Poisson homogeneous spaces that
has been given only recently in \cite{BCST,Lhs} (see also \cite{stef, ferna}). Relying on the explicit 
computation of Poisson cohomology given in \cite{Ginz}, we can conclude that the prequantization cocycle 
is non trivial. We then define a partition of the symplectic groupoid in lagrangian submanifolds whose 
Bohr-Sommerfeld leaves reproduce the Sheu groupoid $\G_S$. This partition, though multiplicative in the sense
of \cite{Hawkins} fails to be a real polarization, even in a generalized sense, due to the singularity of 
the leaf defined by the degenerate point. However it allows to define a topological groupoid on the space of 
Bohr--Sommerfeld leaves which turns out to be exactly the Sheu's groupoid. We then introduce also a complex
polarization on the same groupoid and reconstruct a Hilbert space of states. In both cases the modular function 
preserves the polarization and is quantizable; the corresponding operator is the expected modular operator. 
In the complex case, we discuss a possible convolution product that recovers the convolution algebra
of $\G_S$. It looks more natural to interpret this algebra as an Hilbert algebra, so that the natural outcome of
the construction is the associated Von Neumann algebra. We discuss from this point of view the definition of 
the scalar product, which is a delicate point.

In Section 2 we fix conventions and recall background material. In Section 3 we recall the basic results of 
\cite{WeXu} on the prequantization of a symplectic groupoid and the definition of multiplicative polarization 
given in \cite{Hawkins}. In section 4 we recall basic facts of groupoid $C^*$-algebras following \cite{Renault} 
and discuss in some detail the case of $\O_1$. In Section 5 we collect the semiclassical and quantum results about
the standard Podle\`s sphere, in particular we describe the groupoid $\G_S$ underlying the $C^*$-algebra of the 
quantum space. In Section 6 we describe explicitly the symplectic integration of the Poisson sphere. In Section 7 
we discuss the prequantization and introduce two distinct polarizations: the first one being real and singular 
and which recovers $\G_S$ as groupoid of B-S leaves, the other one being complex which recovers the convolution
algebra on the space of polarized sections.


\section{Notations and conventions}\label{conventions}

We introduce here basic notions and notations of geometric quantization and Lie groupoid theory.

\subsection{Geometric quantization}
We follow notations and conventions of \cite{SN} for hamiltonian dynamics and geometric quantization. 
Let $(M,\Omega)$ be a symplectic manifold. 
For any $f\in C^\infty(M)$ we define the hamiltonian vector field of $f$ as the vector field $\chi_f$ 
such that $\iota_{\chi_f}\Omega = -df$; the Poisson bracket between $f,g\in C^\infty(M)$ is then given by 
$$
\{f,g\}=\langle \Omega,\chi_g\wedge \chi_f\rangle = \chi_g(f)=-\chi_f(g)\, .
$$
Hamiltonian vector fields define a right representation of the Poisson bracket, {\it i.e.} 
$[\chi_f,\chi_g]=-\chi_{\{f,g\}}$.

If the symplectic form is integer, {\it i.e.} it defines a class $\Omega/i\hbar\in H^2(M,\Z)$, there exists  
a {\it prequantization}, {\it i.e.} a line bundle $\Lambda$ over $M$ with connection $\nabla$, whose curvature 
$R_\nabla=-(i/\hbar) \Omega$. If $\Theta$ is a local primitive of $\Omega$, we then write
$\nabla = d -\frac{i}{\hbar} \Theta$.

A {\it polarization} is a lagrangian, involutive distribution $F$ of the complexified tangent bundle $(TM)_\C$. 
It is positive if $i\Omega(X,\bar{X})\geq 0$ for any $X\in \Gamma(F)$. There are two (possibly singular) real 
distribution associated $\D_\C=F\cap\bar{F}$ and $\E_\C=F+\bar{F}$. The polarization is {\it real} if $F=\bar{F}$.

By using the symplectic form we can identify $F^\perp=\Omega(F)$. On $F^\perp$ there is defined the {\it Bott 
connection} $\nabla_X(\xi)=\iota_X d\xi$, for any $X\in\Gamma(F)$. This connection is extended to any square 
root $\sqrt{\det F^\perp}$. The space of polarized sections is 
$$
\H(F)=\{\sigma\otimes \lambda\in \Gamma(\Lambda\otimes\sqrt{\det F})\,| \, \nabla_X(\sigma\otimes\lambda)=0,\, X\in 
\Gamma(F)\}\, .
$$
The natural product between two polarized sections takes value in $\sqrt{\det F^\perp}\otimes
\sqrt{\det\bar{F}^\perp}$, which is naturally isomorphic to $ \sqrt{\det\E_\C^\perp}\otimes
\sqrt{\det \D^\perp_\C}$. The exterior product with $\Omega^k$, where $2k={\rm rk}\E-{\rm rk} \D$, gives the 
isomorphism 
\begin{equation}\label{symplectic_isomorphism}
\Omega^k:~~\det \E^\perp_\C\rightarrow\det {\cal D}^\perp_\C~~~,
\end{equation}
that defines the inner product on $\H(F)$ as an integral on $M/{\cal D}$.
 
Any function $f\in C^\infty(M)$ whose hamiltonian flux preserves the polarization can be quantized. Let
$\underline{b}=\{\chi_i\}$ be a basis of hamiltonian vector fields of $F$; then we have that $[\chi_f,\chi_i]=
\sum_j a_i^j(f) \chi_j$. We then define
\begin{equation}\label{quantization}
\hat{f}(\sigma\otimes\sqrt{\underline{b}})= 
\left[-i\hbar \nabla_{\chi_f}+f-\frac{i\hbar}{2}\sum_i a^i_i(f)
\right]\sigma\otimes\sqrt{\underline{b}}\quad .
\end{equation}

\medskip

\subsection{Symplectic groupoid}
Let $\G=(\G,\G_0,l_\G,r_\G,m_\G,\iota_\G,\epsilon_\G)$ be a Lie groupoid over the space of units $\G_0$, where 
$l_\G, r_\G:{\cal G}\rightarrow \G_0$ are the source and target maps, respectively, $m_\G:\G_{2}\rightarrow\G$ is 
the multiplication, $\iota_\G:\G\rightarrow\G$ is the inversion and $\epsilon_\G:\G_0\rightarrow\G$ is the embedding 
of units. Our conventions are that $(\gamma_1,\gamma_2)\in{\cal G}_{2}$ if $r_\G(\gamma_1)=l_\G(\gamma_2)$. 
We say that $\G$ is source simply connected (ssc) if $l_\G^{-1}(x)$ is connected and simply connected for any 
$x\in \G_0$.

We denote with $\G_k$ the space of strings of $k$-composable elements of $\G$, with the convention that $\G_1=\G$.
The face maps are $d_i:\G_s\rightarrow \G_{s-1}$, $i=0,\ldots s$, defined for $s>1$ as
\begin{equation}\label{face_maps}
d_i(\gamma_1,\ldots\gamma_s) = \left\{ \begin{array}{ll} (\gamma_2,\ldots\gamma_s)& i=0\cr
(\gamma_1,\ldots\gamma_i\gamma_{i+1}\ldots) & 0<i<s\cr
(\gamma_1,\ldots\gamma_{s-1})& i=s\end{array}\right.
\end{equation}
and for $s=1$ as $d_0(\gamma)=l_\G(\gamma)$, $d_1(\gamma)=r_\G(\gamma)$. The simplicial coboundary operator 
$\partial^*:\Omega^k(\G_{s})\rightarrow\Omega^k(\G_{s+1})$ is defined as
$$
\partial^*(\omega) = \sum_{i=0}^s (-)^i d_i^*(\omega)  \;, 
$$
and ${\partial^*}^2=0$.
The cohomology of this complex for $k=0$ is what is called groupoid cohomology.

Let $\G$ and $\G'$ be groupoids, a {\it groupoid morphism} is a pair of maps $F:\G\rightarrow \G '$ and 
$f:\G_0\rightarrow \G_0'$ such that $ l_{\G '}\circ F=f\circ l_\G$, $r_{\G'}\circ F=f\circ r_\G$ and 
$F(m_\G(\gamma_1,\gamma_2))=m_{\G'}(F(\gamma_1), F(\gamma_2))$ for $(\gamma_1,\gamma_2)\in \G_{2}$. 

A {\it symplectic groupoid} is a Lie groupoid, which is equipped with a symplectic form $\Omega_\G$, such 
that the graph of the multiplication is a lagrangian submanifold of $\G \times \G \times \bar\G$, where 
$\bar\G$ means $\G$ with the opposite symplectic structure. An equivalent characterization for a Lie groupoid 
$\cal G$ to be a symplectic groupoid is that the symplectic form be multiplicative, {\it i.e.} 
$\partial^*(\Omega_\G)=0$.
There exists a unique Poisson structure on $\G_0$ such that $l_\G$ and $r_\G$ are Poisson and anti-Poisson 
mappings, respectively. 
A Poisson manifold is said to be {\it integrable} if it is the space of units of a symplectic groupoid.

Let $(M, \pi)$ be a Poisson manifold, where we denote with $\pi$ the bivector defined by the Poisson 
bracket as $\{f,g\}=\pi^{\mu\nu}\partial_\mu(f)\partial_\nu(g)$. As a consequence of the Jacobi identity 
of the Poisson brackets, $d_{LP}(X)=[\pi,X]$ squares to zero; its cohomology is called the Lichnerowicz-Poisson 
cohomology and is denoted with $H_{LP}(M,\pi)$. Two distinguished classes are relevant for what follows. 
The first one is the class $[\pi]\in H^2_{LP}(M,\pi)$ defined by the Poisson tensor itself. Let us assume that 
$M$ is orientable, and let us choose a volume form $V_M$ on $M$. The modular vector field 
$\chi_{V_M}={\rm div}_{V_M} \pi$ is $d_{LP}$-closed; its cohomology class $\chi_{V_M}\in H^1_{LP}(M,\pi)$ 
does not depend on the choice of the volume form and is called the {\it modular class} \cite{Weinstein1997}. 

To a Poisson manifold $(M,\pi)$ there is always associated a topological groupoid $\Sigma(M )$ 
\cite{CatFel,CraFer}. The elements of $\Sigma(M)$ are equivalence classes of {\it cotangent paths} under 
the {\it cotangent homotopy}. A cotangent path is a $C^1$--path $c:I\rightarrow T^*M$ such that 
$$
\pi^\sharp(c(t))=\frac{d}{dt}p(c(t))\,,
$$ 
where $p:T^*M\rightarrow M$ is the cotangent projection and the sharp map $\pi^\sharp:T^*M\rightarrow TM$ 
denotes the contraction with $\pi$. The source and target maps $l\,,r : \Sigma(M ) \rightarrow M$ are given 
by $l([c])=p(c(0))$ and $r([c]) = p(c(1))$, the multiplication in $\Sigma(M )$ is defined by the concatenation 
of cotangent paths.
In general $\Sigma(M)$ is a topological groupoid, but if it is smooth $\Sigma(M)$  carries a symplectic 
structure that makes it the unique (ssc) symplectic groupoid integrating ($M,\pi$). 

Using this description of the symplectic groupoid, it is easy to see that any vector field $\chi$ of $M$ 
that is closed $d_{LP}(\chi)=0$ can be lifted to a groupoid 1--cocycle $F_\chi$ defined as
\begin{equation}\label{lift_vector_field}
F_{\chi}([c]) = \int_0^1 \langle \chi(p(c(t))),c(t)\rangle dt \;.
\end{equation}
In particular the modular vector field $\chi_{V_M}$ is lifted to the {\it modular function} $F_{V_M}$.

\bigskip
\bigskip

\section{Geometric Quantization of Symplectic Groupoids}

We review here the basic facts of the prequantization of a symplectic groupoid as established in \cite{WeXu} 
and give a notion of polarization as in \cite{Hawkins}. 

The prequantization of a symplectic groupoid $\mG (M)$ over $\G_0=M$, seen as a symplectic manifold, is 
given by the pair $(\Lambda, \nabla)$, where $\Lambda$ is a hermitian line bundle over $\mG (M)$ and $\nabla$ 
is a connection whose curvature is proportional to the symplectic form $\Omega_\G$ of $\mG (M)$. Additional 
requirements must be added to take into account the groupoid structure, as showed by P. Xu and A. Weinstein 
in \cite{WeXu}. 

Let $\partial^* \Lambda^*$ denote the (hermitian) line bundle $d_0^* \Lambda^* \otimes d_2^* \Lambda^* \otimes 
d_1^* \Lambda$ over $\G_2(M)$, where $d_i$ denote the face maps defined in (\ref{face_maps}). We give the 
following definition.

\begin{defn} \label{dfn:sympl_gpd_preq}
A prequantization of the symplectic groupoid $\mG(M)$ consists of the triple $(\Lambda, \nabla; \zeta)$
where $(\Lambda, \nabla)$ is a prequantization of $\mG (M)$ as a symplectic manifold and $\zeta$ is
a section of $\partial^* \Lambda^*$ such that:

\begin{itemize}
 \item[$i$)] $\zeta$ has norm one and is multiplicative, {\it i.e.} it satisfies for 
$(\gamma_1,\gamma_2,\gamma_3)\in \G_3(M)$
 \begin{equation} \label{cocycle_property}
  \zeta(\gamma_1, m_\G(\gamma_2, \gamma_3)) \otimes \zeta(\gamma_2,\gamma_3) = \zeta(\gamma_1, \gamma_2) \otimes 
\zeta(m_\G(\gamma_1, \gamma_2),\gamma_3)\; .
 \end{equation} 
 \item[$ii$)] $\zeta$ is covariantly constant, {\it i.e.} if $\Theta_\G$ is a (local) primitive of $\Omega_\G$, 
then $\zeta$ (locally) satisfies
 \begin{equation}
\label{parallel_cocycle}
d\zeta + (\partial^*\Theta_\G) \zeta=0\;. 
\end{equation}
\end{itemize}
\end{defn}

The original formulation of \cite{WeXu} was given in terms of the \emph{Souriau picture}, where the 
prequantization is a principal $\mathbb{S}^1\,$--bundle $p_E: E \rightarrow \mG (M)$ endowed with
a connection $c$ such that $d c = p_E^* \Omega_\G$. The line bundle $\Lambda$ is then recovered 
as the associated bundle $\Lambda = E \times_{\rho_+} \mathbb{C}$ (so that $\Lambda^* =
E \times_{\rho_-} \mathbb{C}$), where $\rho_\pm (e^{i \phi}) \cdot z = e^{\pm i \phi} z$, for all $z
\in \mathbb{C}$. The connection $c$ is flat when restricted to the lagrangian submanifold of units $\epsilon_\G(M) 
\subset \mG(M)$. However, for a generic prequantization, this flat connection can have non trivial holonomy group. 
The existence of the prequantization cocycle $\zeta$ is equivalent to the triviality of the holonomy. Moreover, the 
cocycle defines a groupoid structure on $E$ such that the restriction $c \oplus c \ominus c$ to the graph of the 
multiplication of $E$ is zero.
 
An important consequence of the groupoid prequantization is the existence of an associative multiplication 
between fibres of $\Lambda$ over composable points of $\mG (M)$. In fact, the cocycle property 
(\ref{cocycle_property}) implies that the product $*_\zeta: \Lambda_{\gamma_1} \otimes \Lambda_{\gamma_2} 
\rightarrow \Lambda_{m_\G(\gamma_1,\gamma_2)}$, defined for any $\psi_{\gamma_1} \in \Lambda_{\gamma_1}, 
\phi_{\gamma_2} \in \Lambda_{\gamma_2}$, $(\gamma_1, \gamma_2) \in \G_2(M)$, by

\begin{equation}\label{multiplication}
\psi_{\gamma_1} *_\zeta \phi_{\gamma_2} = \langle \zeta(\gamma_1,\gamma_2), \psi_{\gamma_1} \otimes \phi_{\gamma_2} 
\rangle \, ,
\end{equation}

is associative. We denote as well by $*_\zeta: \Gamma(\Lambda) \otimes \Gamma(\Lambda) \rightarrow
\Gamma(m^*_\G \Lambda)$ the corresponding multiplication between sections of $\Lambda$.

Using results of \cite{WeXu}, we can state the existence and uniqueness of the prequantization of $\mG(M)$.

\begin{thm} \label{thm:UniPre}
Let $\mG(M)$ be (ssc) and $l_\G$-locally trivial. If $\mG(M)$ is prequantizable as a symplectic manifold, then 
there exists a unique prequantization without holonomy on $\epsilon_\G(M)$.
\end{thm}

The second step in geometric quantization is the choice of a polarization. Also at this step a compatibility 
with the groupoid structure can be required. 

We use the following definition given in \cite{Hawkins} in an attempt to construct a convolution product 
between polarized sections.

\begin{defn} \label{dfn:mulpol}

For any distribution $\mP \subset T_\C \mG(M)$, let us denote $\mP_2 = (\mP \times  \mP) \cap T_\C \G_2(M)$.
Then $\mP$ is called multiplicative if, for every $(\gamma_1, \gamma_2) \in \G_2(M)$,
$$
  m_{\G*} \, (\mP_2)_{(\gamma_1,\gamma_2)} = \mP_{m_\G (\gamma_1, \gamma_2)}.
$$
\end{defn}

Let us now define, for a polarized section $\psi$, an adjoint map by $\psi^\dagger:= \iota_\G^* 
(\overline{\psi})$. Then, requiring that the symplectic potential of $\nabla$ is antinvariant with 
respect to $\iota_\G$, it can be demonstrated that $\psi^\dagger$ is still polarized if the distribution 
is {\it hermitian}, \textit{i.e.} $ \iota_{\G*} \,  \mP = \overline{\mP}$.
We are thus led to:

\begin{defn} \label{def:SymplGpoidPolar}
A polarization $\mathcal{P }\subset T_\C \mG(M)$ of $\mG(M)$ as symplectic manifold is called
a polarization of $\mG(M)$ as a symplectic groupoid if it is multiplicative and hermitian. 
\end{defn}

\medskip
\begin{rem}{\rm The basic model of a real multiplicative polarization is given as the kernel 
of a groupoid fibration. This is what happens for instance in  \cite{We934}.}
\end{rem}

\bigskip
\bigskip

\section{Groupoid C$^*$-algebras}

We informally describe some basic facts of the construction of groupoid $C^*$-algebras, following 
\cite{Renault}. Let $\G$ be a topological groupoid and let $C_c(\G)$ denote the space of continuous 
functions with compact support. A {\it left Haar system} for $\G$ is a family of measures 
$\{\lambda^x, x\in \G_0\}$ on $\G$ such that
\begin{itemize}
\item[$i$)] the support of $\lambda^x$ is $\G^x=l^{-1}_\G(x)$;
\item[$ii$)] for any $f\in C_c(\G)$, and $x\in \G_0$, $\lambda(f)(x)=\int_\G f d\lambda^x$, defines 
$\lambda(f)\in C_c(\G_0)$;
\item[$iii$)] for any $\gamma\in\G$ and $f\in C_c(\G)$, $\int_\G f(\gamma\gamma')d
\lambda^{r_\G(\gamma)}(\gamma')=$ $\int_\G f(\gamma') 
d\lambda^{l_\G(\gamma)}(\gamma')$.
\end{itemize}

The composition of $\lambda^x$ with the inverse map will be denoted as $\lambda_x$; this family 
defines a {\it right Haar system} $\lambda^{-1}$. 
Any measure $\mu$ on the space of units $\G_0$ induces measures $\nu,\nu^{-1}$ on the whole  $\G$ 
through 
$$
\int_\G f d\nu = \int \lambda(f)d\mu \; ,\qquad\int_\G f d\nu^{-1}=\int \lambda^{-1}(f)d\mu\;.
$$
The measure $\mu$ is said to be {\it quasi-invariant} if $\nu$ and $\nu^{-1}$ are equivalent measures; 
in this case the Radon-Nikodym derivative $D=d\nu/d\nu^{-1}$ is called the {\it modular function} of $\mu$. 
The function $\log D\in Z^1(\G,\R)$ turns out to be a groupoid one cocycle with values in $\R$ and its 
cohomology class depends only on the equivalence class of $\mu$.

The notion of quasi-invariant measure can be extended as follows: let $c\in Z^1(\G,\R)$ be a fixed cocycle 
and $\beta\in[-\infty,\infty]$. Consider the set:
$$
{\rm Min}(c)=\{x\in\G_0\, :\, c(\G_x)\subset [0,\infty)\, \}\; .
$$
A measure $\mu$ on $\G_0$ satisfies the $(c,\beta)$--{\it KMS condition} if
\begin{itemize}
\item[$i$)] when $\beta$ is finite, $\mu$ is quasi invariant and the modular function is $e^{-\beta c}$.
\item[$ii$)] when $\beta=\pm\infty$, the support of $\mu$ is contained in ${\rm Min}(\pm c)$.
\end{itemize}
A normalized $(c,\infty)$--KMS measure is called a \emph{ground state} for $c$.

Let $\zeta\in Z^2(\G,\T)$ be a continuous two-cocycle. For any $f,g \in C_c(\G)$ let us define the convolution 
and the involution as
\begin{eqnarray}\label{groupoid_convolution}
(f*g) (\gamma) &:=& \int f(\gamma\gamma')g(\gamma'{}^{-1})\zeta(\gamma\gamma',\gamma'{}^{-1})\ d
\lambda^{r_\G(\gamma)}(\gamma')\;,\cr
f^*(\gamma) &:=& \overline{f(\gamma^{-1})}\ \overline{\zeta(\gamma,\gamma^{-1})} \;\;
\end{eqnarray}
The space $C_c(\G)$ equipped with these operations defines the $*$-algebra $C_c(\G,\zeta)$. 
The $C^*$-norm is defined as $||f||=\sup_L ||L(f)||$ over all bounded representations $L$ (we skip details, 
see \cite{Renault}). The completion of $C_c(\G,\zeta)$ with respect to this norm defines the $\zeta$-twisted 
convolution $C^*$-algebra $C^*(\G,\zeta)$ of the groupoid $\G$. We omit $\zeta$ in the notation when the cocycle 
is the trivial one.

Let now $c\in Z^1(\G,\R)$ be a one cocycle with values in $\R$. It defines an algebra automorphism 
$$
A_c:\R\rightarrow {\rm Aut}(C_c(\G,\zeta))~~,\quad (A_c(t)f)(\gamma)= e^{it c(\gamma)}f(\gamma)~~~.
$$
A measure $\mu$ on the space of units $\G_0$ defines the weight $\phi_\mu$ on $C^*(\G,\zeta)$ as 
\begin{equation}
\label{KMS_weight} 
\phi_\mu(f)=\int_{\G_0} f d\mu~~~.
\end{equation} 
We recall that given the automorphism $A_c$ the weight $\phi_\mu$ satisfies the KMS condition at 
$0\leq\beta<\infty$ if, for any $f,g\in C_c(\G,\zeta)$, we have 
$$
\phi_\mu(f * A_c(i\beta)(g)) = \phi_\mu(g * f)  ~~~~ ;
$$
if $\beta=\infty$, $\phi_\mu$ satisfies the KMS condition at $\infty$ if $|\phi_\mu(f * A_c(z)(g))|\leq ||f|| 
\ ||g||$, for any $z$ such that ${\rm Im}z>0$.
It is proven that $\mu$ satisfies the $(c,\beta)$-KMS condition for $\beta\in[0,\infty]$ if and only if $\phi_\mu$ 
satisfies the KMS condition for the automorphism $A_c$ at $\beta$. 

The GNS representation generated by the KMS weight $\phi_\mu$ is obtained as convolution action of $C_c(\G,\zeta)$ 
on the Hilbert space $L^2(\G,\nu^{-1})$. If $\beta$ is finite (and so $\mu$ is quasi invariant and $D=e^{-\beta c}$ 
is the modular function) then $C_c(\G,\zeta)$, equipped with the inner product of $L^2(\G,\nu^{-1})$, is a left 
Hilbert algebra. We recall the definition.
 
\begin{defn} 
A $*$-algebra ${\cal A}$, equipped with a inner product, is a left Hilbert algebra if the left regular 
representation is bounded and involutive, if $S(f)=f^*, f \in {\cal A}$, is preclosed, and if ${\cal A}^2$ is 
dense in ${\cal A}$.  
\end{defn}

The polar decomposition of $S=J\Delta^{1/2}$ defines the modular conjugation $J$ and the modular operator 
$\Delta=S^\dagger S$. On $f\in C_c(\G,\zeta)$ we have
$$
(J f)(\gamma)= D^{1/2}(\gamma)f^*(\gamma)~,~~~~~
(\Delta f)(\gamma) = D(\gamma) f(\gamma)\;.
$$
In case $\beta=\infty$ the KMS condition means that the hamiltonian $c$ is positive on ${\rm supp}(\nu)$ 
and the associated state $\phi_\mu$ is called a {\it ground state}.

\subsection{The groupoid $\O_1$}\label{ToplitzGrpd}

Let us discuss here an example of the previous construction that will be relevant in the following sections. 
Let $\Z$ act on $\overline{\Z}=\Z\cup \{\infty\}$ by translation and leaving $\infty$ fixed, and let 
$\overline{\Z}\times\Z$ be the action groupoid. Here $\Z$ is considered with the discrete topology and 
$\overline\Z$ is the compactification to $+\infty$, therefore the action groupoid is a locally compact
Hausdorff topological groupoid. The $n=1$ Cuntz groupoid (see \cite{Renault}) is the restriction 
$\O_1=(\overline{\Z}\times\Z)|_{\overline{\N}}$ of the action groupoid to $\overline{\N}\subset\overline{\Z}$. 
The Haar system is given by the discrete measure on the $l$-fibre and the convolution algebra $C_c(\O_1)$ 
is defined by
\begin{equation}\label{cuntz_grpd_relations}
e_{m,n}*e_{p,q} = \delta_{m+n,p} e_{m,n+q} \;,~~~~~~~ e_{m,n}^* = e_{m+n,-n} \;,
\end{equation}
where $e_{m,n}(p,q)=\delta_{mp}\delta_{nq}$. The groupoid $C^*$-algebra $C^*({\cal O}_1)$ can be described 
as the $C^*$-algebra generated by the shift operator ${\cal S}$ on $\ell^2(\N)$, the explicit isomorphism 
being ${\cal S}=\sum_{m\geq 0} e_{1+m,-1}$, with the action $e_{m,n}|p\rangle = \delta_{m+n,p}|m\rangle$.

Let us consider the cocycle $c_1(m,n)=n$; the corresponding automorphism is
\begin{equation}
\label{modular_automorphism}
A_{c_1}(t)e_{m,n}=e^{itn} e_{m,n}\;.
\end{equation}
In \cite{Renault} it is shown that there exists a unique $(c_1,\beta)$-KMS probability measure on 
$\epsilon(\O_1)=\overline{\N}$ for $\beta\leq 0$: for $\beta=-\infty$ it reads $\mu_{-\infty}(0)=1$, 
for $\beta=-\hbar<0$ it is $\mu_\hbar(n)=e^{-n\hbar }(1-e^{-\hbar})$ and for $\beta=0$ it is $\mu_0(\infty)=1$.

Let us analyse the corresponding GNS constructions. 

Let us consider first the measure $\mu_\hbar$; the GNS state (\ref{KMS_weight}) is given by
\begin{equation}\label{measure_h}
\phi_{\mu_\hbar}(f)=\sum_m f(m,0)\mu_\hbar(m)~~~. 
\end{equation}
The GNS Hilbert space $\H_{\mu_\hbar}$ has a basis given by $\{|m,n\rangle \equiv e_{m,n-m}\}$ with 
$\langle m,n|m',n'\rangle=\phi_{\mu_\hbar}(e_{m,n-m}^**e_{m',n'-m'})=\delta_{mm'}\delta_{nn'} \mu_\hbar(n)$ 
and action
$$
e_{m,n-m}|p,q\rangle = \delta_{np}|m,q\rangle\;.
$$
The cyclic vector is $|\xi\rangle = \sum_{n\geq0}|n,n\rangle$, the modular operator and conjugation respectively 
are $\Delta|m,n\rangle = e^{-\hbar(m-n)}|m,n\rangle$ and $J|m,n\rangle = e^{-\hbar(m-n)/2}|n,m\rangle$. 
 
The GNS state corresponding to the ground state $\mu_{-\infty}$ is 
$\phi_{\mu_{-\infty}}(f)=f(0,0)$, and the Hilbert space is ${\cal H}_{\mu_{-\infty}}=\ell^2(\N)$ with action 
$$e_{m,n}|p\rangle = \delta_{m+n,p}| m \rangle~.$$

Finally the GNS state corresponding to $\mu_0$ is $\phi_{\mu_0}(f)=f(\infty,0)$ and the GNS Hilbert space 
is ${\cal H}_{\mu_0}=L^2(\sphere^1)$ with action given by $f F=\sigma(f) F$, where $F\in L^2(\sphere^1)$ and 
$\sigma:C^*(\O_1)\rightarrow C(\sphere^1)$ is defined as
\begin{equation}\label{evaluation}
\sigma(f)(x) = \sum_m f(\infty,m) x^m \;.
\end{equation}
One can easily check that $\sigma({\cal S})(x)=x^{-1}$. Finally, if we denote with ${\cal K}$ the 
$C^*$-algebra of compact operators, it is shown in \cite{Co} that the following sequence is exact
\begin{equation}\label{exact_sequence}
0\rightarrow {\cal K} \rightarrow C^*(\O_1) \overset\sigma\rightarrow C(\sphere^1) \rightarrow 0\;.
\end{equation}
The geometrical picture of (\ref{exact_sequence}) consists of looking to $C^*(\O_1)$ as the quantization 
of a closed disk whose Poisson structure is symplectic in the interior and zero on the boundary.

\bigskip
\bigskip

\section{The Podle\`s spheres}\label{podles_sphere}
The ``quantum spheres''are a family parametrized by $c\in\R$ of unital $*$-algebras ${\cal A}(\qsfera)$ 
introduced in \cite{Podles}; the deformation parameter $q$ is real. They are generated by $\tau=\tau^*$ and 
$\alpha$ with relations
$$ q^2\alpha^*\alpha=\tau(1-\tau)+c I,~~~ 
q^2 \alpha \alpha^*=q^2\tau (1-q^2\tau)+c I,~~~  \alpha\tau=q^2\tau\alpha \ .
$$
We are interested in the case $c=0$. There are only two irreducible representations, the counit $\varepsilon$ 
and the infinite dimensional $\rho$. The counit is $\varepsilon(\tau)=1,~\varepsilon(\alpha)=0$. With respect 
to an orthonormal basis $\{\psi_n\}_{n=0}^\infty$, $\rho$ is defined as
($0<q<1$):
\begin{equation}\label{unitary_rep_podles_sphere}
\rho(\tau)\,\psi_n=q^{2n}\,\psi_n~,~~~~~~\rho(\alpha)\,\psi_n=q^{n-1}(1-q^{2n})^{1/2}\,\psi_{n-1}\;.
\end{equation}

The operators $\rho(\alpha)$ and $\rho(\tau)$ are compact and the associated $C^*$-algebra
 $C({\mathbb S}^2_{q,0})$ is isomorphic to the minimal unitization of the compact operators.
 
In \cite{Sheu1} it is proved that the C*-algebras of quantum spheres can be realized as concrete 
groupoid C*-algebras. In terms of the generators $e_{m,n}$ of $C^*(\O_1)$, discussed in the previous
 section, the generators of ${\cal A}({\mathbb S}^2_{q,0})$ are:
\begin{equation}\label{qsfera_inside_cuntz}
\tau=\sum_{m\geq 0} q^{2m} e_{m,0}\,,\qquad \alpha=\sum_{m\geq 0} q^m(1-q^{2(m+1)})^{1/2} e_{m,1}\,.
\end{equation}
They generate the $C^*$-subalgebra $C^*(\sphere^2_{q,0})=\sigma^{-1}(\C)$, where $\sigma$ is defined in 
(\ref{evaluation}); since $C^*(\O_1)$ is the quantization of the disk whose Poisson structure is symplectic 
in the interior and zero on the boundary, the $C^*$-subalgebra $\sigma^{-1}(\C)$ corresponds to identify the 
boundary to a point. In the groupoid picture this corresponds to restrict to the subgroupoid 
\begin{equation}\label{sheu_groupoid} \G_S=\{(m,n)\in\O_1\ | \ m=\infty\implies n=0\}\subset\O_1. 
\end{equation}
The Sheu groupoid $\G_S$ inherits from $\O_1$ the Haar system, so that $C^*(\sphere^2_{q,0})=C^*(\G_S)$, and 
the cocycle $c_1$ whose automorphism (\ref{modular_automorphism}) reads $A_{c_1}(t)\tau=\tau$ and $A_{c_1}(t)
\alpha= e^{it}\alpha$. Finally the state $\phi_{\mu_\hbar}$ given in (\ref{measure_h}) is the usual Haar measure 
on $\sphere^2_{q,0}$ when we identify $q=e^{-\hbar/2}$; as a KMS-state, $\phi_{\mu_\hbar}$ satisfies
$$
\phi_{\mu_\hbar}(f *A_{c_1}(-i\hbar)g) = \phi_{\mu_\hbar}(g*f)\;.
$$ 
Let us recall that $A_{c_1}(-i\hbar)f=D(f)$ with $D$ being the modular operator discussed in Section 
\ref{ToplitzGrpd}.

The Poisson structure ($\sphere^2,\pi$) is the semiclassical limit of  ${\cal A}({\mathbb S}^2_{q,0})$. 
In terms of the complex generator $\alpha$ and real $\tau$ satisfying $|\alpha|^2=\tau(1-\tau)$ the bracket 
is 
$$
\{\alpha,\tau\}=i \alpha \tau,\qquad \{\alpha,\bar \alpha\}= i\tau (1-2 \tau).
$$
The symplectic foliation is made of a singular point $\{N\}$ defined by $\tau=0$ and of a two dimensional 
symplectic leaf on the complement. On $\sphere^2\setminus\{N\}$, the Poisson tensor reads in terms of the 
complex coordinate $z=\alpha/\tau$ as
$$
\pi = -i(1+|z|^2)\partial_z\wedge\partial_{\bar z}\,;
$$
on the chart $\sphere^2\setminus\{S\}$, where $\{S\}$ is defined by $\tau=1$, we have in terms of $w=1/z=\bar 
\alpha/(1-\tau)$ 
$$
\pi =-i |w|^2 (1+|w|^2)\partial_w\wedge\partial_{\bar w}\,.
$$
The Poisson cohomology of $(\sphere^2,\pi)$ has been computed in Proposition 2.18 of \cite{Ginz}. The result 
is
\begin{equation}
 \label{Poisson_cohomology_podles_sphere}
 H^{(0)}_{LP} = \R ~,~~~ H^{(1)}_{LP} = \R[\chi_{V_{\sphere^2}}]~,~~~ H^{(2)}_{LP} = \R[\pi] \oplus 
\R[i\partial_z\wedge\partial_{\bar z}]~,
\end{equation}
where $\chi_{V_{\sphere^2}}$ is the modular vector field that we will describe in (\ref{modular_vector_field}).

\bigskip
\bigskip

\section{The symplectic groupoid of ($\sphere^2,\pi$)}
We describe here the explicit integration of ($\sphere^2,\pi$). Since ($\sphere^2,\pi$) is a Poisson 
homogeneous space, {\it i.e.} $\sphere^2=U(1)\backslash SU(2)$, with $SU(2)$ equipped with the Poisson Lie 
group structure, we can use the general construction given in \cite{BCST}.

The general framework for Poisson homogeneous spaces is the following, see \cite{LuTh} for a detailed 
account of Poisson reduction. Let $G$ be a simply connected Poisson Lie group and let $G^*$ be the dual 
Poisson-Lie group of $G$, that is simply connected by definition. We assume that $G$ is complete 
and we denote the left and right {\it dressing action} of $G$ on $G^*$ as ${}^g\gamma$ and $\gamma^g$, 
respectively, for $g\in G$ and $\gamma\in G^*$. 
Analogously we denote with ${}^\gamma g$ and $g^\gamma$ the dressing actions of $G^*$ on $G$. We recall 
that the double group is $G\times G^*$ as a manifold and that the dressing actions are defined by the following 
relations in $D$: $g\gamma={}^g\gamma g^\gamma$ and $\gamma g= {}^\gamma g \gamma^g$, for $g\in G$ and 
$\gamma\in G^*$.

Let $H\subset G$ be a closed coisotropic subgroup; then on $H\backslash G$ there exists a canonical Poisson 
structure $\pi_{H\backslash G}$ such that the quotient map is Poisson. We denote the quotient map 
$G\rightarrow H\backslash G$ as $g\rightarrow[g]$. If we denote $\g=Lie(G)$ and $\h=Lie(H)$ the annihilator 
$\h^\perp\subset\g^*$ is a subalgebra as a consequence of coisotropy; let $H^\perp\subset G^*$ be the coisotropic 
subgroup integrating $\h^\perp$.

The general result established in \cite{BCST} describes a symplectic groupoid of the Poisson homogeneous 
space $(H\backslash G,\pi_{H\backslash G})$ as
$$
{\cal G}(H\backslash G,\pi_{H\backslash G})=\{([g],{}^{g^{-1}}\sigma)\in H\,
\backslash G \times G^*\,:\, 
\sigma\in H^\perp\} \;.
$$
If we consider a (local) section $s_U:U\subset H\backslash G\to G$, we define the 
trivialization
$U\times H^\perp\rightarrow{\cal G}(H\backslash G,\pi_{H\backslash G})\big|_U$ as 
$$
( [g],\sigma)\mapsto ([g],{}^{s_U[g]^{-1}}\sigma).
$$

If $H$ is a Poisson Lie subgroup and not just coisotropic, then the exponential map 
$\exp:\h^\perp\rightarrow H^\perp$ intertwines the adjoint with the dressing action. If moreover 
$\exp$ is a diffeomorphism, {\it i.e.} $H^\perp$ is of exponential type, then the symplectic groupoid 
$\G(H\backslash G,\pi_{H\backslash G})$ is diffeomorphic to $T^*(H\backslash G)$. 

We are going to apply this description to the standard Podle\`s sphere $\sphere^2=U(1)\backslash SU(2)$, 
where $G=SU(2)$ is equipped with the Poisson structure corresponding to the factorization $D=SL(2,\C)=SU(2)
\times SB(2,\C)$, where
$$
SB(2,\C)=\{\left(\begin{array}{cc}v & n \cr 0 & v^{-1}\end{array}\right)~,~~~ v>0, n\in\C\}\;.
$$
Moreover the subgroup $H$ is the diagonal $U(1)$ and $H^\perp$ is the subgroup of $SB(2,\C)$ corresponding 
to $v=1$ in the above parametrization. Here $H^\perp$ is of exponential type with $H^\perp=\exp\h^\perp=1+
\h^\perp\simeq\C$. 
We call $N$, the north pole, the zero symplectic leaf of ($\sphere^2,\pi$) and $S$ the opposite one:
$$
N=\left[\left(\begin{array}{cc}1&0\\ 0&1\end{array}\right)\right]\,,\qquad S=\left[\left(\begin{array}{cc}0&1\\ 
-1&0\end{array}\right)\right] ~~~~~~~.
$$

Let us choose on $\sphere^2\setminus\{S\}$ ( the {\it singular chart}) the complex coordinate
$$
w\left[\left(\begin{array}{cc}a&b\\ -b^*&a^*\end{array}\right)\right] = \frac{b}{a}
$$
and the section
$$
g_S(w)=\frac 1{\sqrt{1+|w|^2}}\left(\begin{array}{cc}1&w\\ -\bar w&1\end{array}\right)\,,
$$
so that we have the following trivialization:
$$
(w,\sigma_S)\mapsto \left(w, {}^{g_S(w)^{-1}}\left(\begin{array}{cc}1&\sigma_S\\ 0&1\end{array}\right)\right)\,.
$$
On $\sphere^2\setminus\{N\}$ (the {\it symplectic chart}) we choose the complex coordinate
$$
z\left[\left(\begin{array}{cc}a&b\\ -b^*&a^*\end{array}\right)\right] = \frac{a}{b}
$$
and the section
$$
g_N(z)=\frac 1{\sqrt{1+|z|^2}}\left(\begin{array}{cc}z&1\\ -1&\bar z \end{array}\right)\,,
$$
that gives the trivialization:
$$
(z,\sigma_N)\mapsto \left(z, {}^{g_N(z)^{-1}}\left(\begin{array}{cc}1&\sigma_N\\ 0&1\end{array}\right)\right)\,.
$$
We observe that
$$
g_S=\frac 1{|w|}\left(\begin{array}{cc}w&0 \\ 0&\bar w\end{array}\right)g_N=t_{SN}g_N\,.
$$
From the condition that the two charts coincide on the overlap we get:
$$
{}^{g_N^{-1}}\left(\begin{array}{cc}1&\sigma_N\\ 0&1\end{array}\right)=
{}^{g_S^{-1}t_{SN}}\left(\begin{array}{cc}1&\sigma_N\\ 0&1\end{array}\right)
={}^{g_S^{-1}}\left(\begin{array}{cc}1&\frac {w^2}{|w|^2}\sigma_N\\ 0&1\end{array}\right)\,,
$$
then
$$
\sigma_S=\frac{w^2}{|w|^2}\sigma_N\,.
$$
Defining
$$
p_N=\frac{\bar\sigma_N}{1+|z|^2},\qquad p_S=-\frac{\bar\sigma_S}{1+|w|^2}\ ,
$$
we have that $p_N$ and $p_S$ transform as the coordinates of the cotangent bundle, {\it i.e.}
$p_N=-w^2p_S$.

This makes explicit the diffeomorphism between ${\cal G}(\sphere^2,\pi_{H\backslash G})$ and $T^*\sphere^2$. 
Let us now discuss the groupoid structure in the two charts. It is clear what it means for $l$ defined 
as $l([g],\gamma)=[g]$; we have to compute the formulas for $r$ defined as $r([g], \gamma)= [g^\gamma]$. 
We observe that defining $\sigma={}^g\gamma \in H^\perp\,,$ we have that $g\gamma={}^g\gamma g^\gamma=\sigma 
g^\gamma={}^\sigma(g^\gamma)\,\sigma^{g^\gamma}$ then $g={}^\sigma(g^\gamma)$ or
${}^{\sigma^{-1}}g=g^\gamma$. In the double quotient $G\times G^* \rightarrow H\backslash G\times G^*/G^*
\simeq H\backslash G$ the element  $\sigma^{-1} g={}^{\sigma^{-1}}g(\sigma^{-1})^g\in G\times G^*
\rightarrow [{}^{\sigma^{-1}}g]=[g ^\gamma]$. Then for computing $[g^\gamma]$ in terms of the 
trivialization coordinates we can compute the coordinate of the class of $\sigma^{-1}g$ in the double quotient 
$G\times G^*\rightarrow H\backslash G$. We have that
$$
w\left[\left(\begin{array}{cc}A&B\\ C&D\end{array}\right)\right] = -\bar C/A\;~~~~~ 
\left(\begin{array}{cc}A&B\\ C&D\end{array}\right)\in SL(2,\C).
$$

An explicit computation with $g=g_S$ and $\sigma=\left(\begin{array}{cc}1&\sigma_S\\ 0&1\end{array}\right)$ 
gives in the singular chart,
$$
l(w,\sigma_S)= w~,~~~~~ r(w,\sigma_S)=[g_{\,S}(w)^\gamma]=\frac w{1+\sigma_S\bar w}
$$
and, in the symplectic chart with $g=g_N$ and $\sigma=\left(\begin{array}{cc}1&\sigma_N\\ 0&1\end{array}\right)$,
\begin{equation}\label{target_grpd_pdls_sphere}
l(z,\sigma_N)=z~,~~~~~~r(z,\sigma_N)=[g_{\,N}(z)^\gamma]=z+\sigma_N  \;.
\end{equation}

The chart without the singularity is a simply connected symplectic manifold that integrates to the pair 
groupoid; there exists a symplectic groupoid morphism $\phi:{\cal G}(\sphere^2)\big|_{\sphere^2\setminus\{N\}}
\to\C\times\C$ which is given by
\begin{equation}\label{diffeo1}
\phi(z,p_N)=(z,z+\sigma_N(p_N,z))=\left(z,z+(1+|z|^2)\bar p_N\right)\,,
\end{equation}
as a consequence of commutativity with source and target maps. We use it to read the multiplication 
rule in ${\cal G}(\sphere^2)$; for $z'=z+\sigma_N=z+(1+|z|^2)\bar{p}_N$ we have
\begin{equation}
\label{multiplication_grpd_pdls_sphere}(z,p_N)(z',p_N')=
(z, p_N + \frac{1+|z'|^2}{1+|z|^2}p_N')\;.
\end{equation}

Analogously, the symplectic form in this chart can be read from the pair groupoid as
\begin{equation}\label{symplectic_form_grpd_pdls_sphere}
\Omega = \frac {1}{i(1+|z|^2)}dz\wedge d\bar z-\frac {1}{i(1+|z+\sigma_N|^2)}(dz+d\sigma_N)
\wedge(d\bar z+d\bar \sigma_N)\,.
\end{equation}
We summarize all the discussion in the following Proposition.

\smallskip
\begin{prop}
\label{symplectic_grpd_pdls_sphere}
The (ssc) symplectic groupoid $\G(\sphere^2,\pi)$ is diffeomorphic to $T^*\sphere^2$, equipped with 
the source and target maps (\ref{target_grpd_pdls_sphere}), the multiplication 
(\ref{multiplication_grpd_pdls_sphere}) and the symplectic form (\ref{symplectic_form_grpd_pdls_sphere}).
\end{prop}

\medskip

Let $V_{\sphere^2}=idz d\bar{z}/(1+|z|^2)^2$ be the round volume form on $\sphere^2$. The 
{\it modular vector field} $\chi_{V_{\sphere^2}}$ reads
\begin{equation}\label{modular_vector_field}
\chi_{V_{\sphere^2}}= i(z\partial_z -\bar{z}\partial_{\bar{z}})\;
\end{equation}
and defines a non trivial class in Poisson cohomology, as seen in Section \ref{podles_sphere}.

The integrated counterpart is the {\it modular function} $c_{V_{\sphere^2}}\in C^\infty(T^*\sphere^2)$ 
defined as
$$
c_{V_{\sphere^2}}(x,y) = \int_0^1 dt \langle \eta(t), \chi_{V_{\sphere^2}}(\gamma(t))\rangle\;,
$$
for any cotangent path $(\gamma,\eta)$ such that $\gamma(0)=x$, $\gamma(1)=y$. It is a smooth function and 
defines a real valued groupoid cocycle, {\it i.e.} $\partial^*c_{V_{\sphere^2}}=0$.  After an easy computation 
we get 
\begin{eqnarray}
\label{symplectic_modular_function}
c_{V_{\sphere^2}}(x,y) = \log D_{V_{\sphere^2}}&=& \log \frac{1+|y|^2}{1+|x|^2} = \log \frac{1+|z+(1+|z|^2)
\bar{p}_N|^2}{1+|z|^2} \cr 
&=& \log\frac{|w|^2 + |1-(1+|w|^2)\bar{w}\bar{p}_S|^2}{1+|w|^2}\;,
\end{eqnarray}
where $x=z$ and $y=z+\sigma_N$. When restricted to the symplectic chart, it is trivial, indeed 
$c_{V_{\sphere^2}}=\partial^*\phi$ with $\phi(x)=\log(1+|x|^2)$. Since under the Van Est map it 
corresponds to the modular class, as a cocycle of the whole $\G(\sphere^2,\pi)$, it is not trivial.

In terms of the groupoid structures the interpretation of $c_{V_{\sphere^2}}$ goes as follows.
Up to a scalar multiplication, any Haar system on $\G(\sphere^2,\pi)$ is written as 
\begin{equation}\label{haar_system_podles}
\lambda^x= -i\Lambda(y)dy d\bar{y},\ (|x|< \infty);~~~~~~~~ \lambda^{\infty}=-i dp_S d\bar{p}_S ~,
\end{equation}
for any positive $\Lambda$ such that  $\lim_{y\rightarrow\infty}\Lambda(y)=1$. If we choose on $\sphere^2$ 
the round volume form $V_{\sphere^2}$, the induced volume form on $\G(\sphere^2,\pi)$ is $\nu(V_{\sphere^2})
 = \Lambda D_{V_{\sphere^2}}\Omega^2$, so that
$$
\nu = \nu^{-1} D_{V_{\sphere^2}}^{2} e^{-\partial^*\log\Lambda} \;,
$$
and the modular cocycle is $2 c_{V_{\sphere^2}}- \partial^* \log\Lambda$.

\begin{rem}
{\rm ($i$) If we restrict $\G(\sphere^2,\pi)$ to the singular chart we get a symplectic groupoid integrating 
the Poisson structure on $\C\simeq\R^2$ with quadratic tensor $\pi=i |w|^2(1+|w|^2)\partial_w\wedge 
\partial_{\bar w}$. The restriction on the symplectic $\C^*=\C\setminus\{0\}$ is the pair groupoid 
$\C^*\times\C^*$, so that the source fibre is not simply connected. The (ssc) groupoid integrating 
this Poisson structure was obtained in \cite{CDH}.

($ii$) Since $\G(\sphere^2,\pi)$ is a Lie groupoid the modular cocycle can be computed without introducing 
the Haar system, as explained in \cite{Weinstein1997,ELW}. Both computations agree and define a class in 
groupoid cohomology that is twice the one defined in (\ref{symplectic_modular_function}), in the same way 
the algebroid modular class of $T^*M$ is twice the Poisson modular class (see \cite{ELW}).}
\end{rem}

\bigskip
\bigskip

\section{Geometric quantization of $\G(\sphere^2,\pi)$}

Let us summarize what can be said about  the prequantization of $\G(\sphere^2,\pi)$. The symplectic form 
$\Omega$ in (\ref{symplectic_form_grpd_pdls_sphere}) is exact; this can be easily seen by looking that its 
restriction on the source fibres is exact and by applying Corollary 5.3 of \cite{Crainic}. The prequantization
 line bundle is therefore the trivial one $T^*\sphere^2\times\C$ with connection equal to $(i/\hbar)\Theta$, 
for some choice of  $\Theta$ such that $d\Theta=\Omega$. The prequantization cocycle $\zeta_\Theta:
\G_2(\sphere^2,\pi)\rightarrow\sphere^1$ satisfies (\ref{parallel_cocycle}) so that it  depends on this choice.
 In particular, we can ask if there exists  a multiplicative primitive of $\Omega$, {\it i.e.} a primitive $\Theta$
 such that $\partial^*\Theta=0$. This would imply $\zeta_\Theta=1$.
By applying Theorem 4.2 in \cite{Crainic} and the non triviality of the class of the Poisson tensor $\pi$ in 
Poisson cohomology, that we mentioned in Section \ref{podles_sphere}, we conclude that the answer is negative.

Let us now come to the choice of a symplectic groupoid polarization. Our guiding principle here will be that 
since we would like to be able to quantize the modular function $F_V$, our lagrangian polarizations should be 
$F_V$--invariant. We will propose two different choices. The first one is a real singular polarization 
(topological constraints here does not allow a non singular one) and then a complex one.

\medskip
\subsection{A real polarization}

In this case our strategy will be to choose first a polarization on the symplectic chart and then try to 
extend it to the whole groupoid ${\cal G}(\sphere^2,\pi)$.

The map 
$\phi:{\cal G}(\sphere^2,\pi)\big|_{\sphere^2\setminus\{N\}}\to
\C\times\C$ given in (\ref{diffeo1}) identifies two variables $(x,y)$ as:
$$
x=z=\frac{1}{w}\;,\;\;\;\;\; y = z + (1+|z|^2)\bar{p}_N = \frac{1}{w}[1-(1+|w|^2)\bar{w}\bar{p_S}]\;,
$$
with inverse given in the symplectic chart by
$$
z= x~,~~~ p_N= \frac{\bar y-\bar x}{1+|x|^2}~,
$$
and in the singular chart by
$$
w= 1/x~,~~~ p_S= -x^2\frac{\bar y-\bar x}{1+|x|^2} \\\\\\ .
$$
Let us consider the harmonic oscillator polarization of $\C\times\C$ given by the contour levels of 
$F_+(x)=|x|^2$ and $F_-(y)=|y|^2$. This is a real multiplicative polarization of ${\cal G}(\sphere^2,\pi)$
restricted to the chart. To be precise, we have already introduced  harmonic oscillator singularities, though 
they can be considered to be not relevant from the point of view of quantization (see the discussion about 
elliptic singularities in \cite{Ham}). If we add the missing lagrangian leaf $T^*_N\sphere^2$ we get a partition 
of $T^*\sphere^2$ in disjoint Lagrangian submanifolds.

Unfortunately, this partition fails to define a continuous distribution (even in a generalized sense) of the whole 
$\G(\sphere^2,\pi)$. In order to understand the topology around the lagrangian leaf $\infty\equiv T^*_N\sphere^2$, 
it is better to describe the space of leaves in the following way. Let $\R$ act on $\bar{\R}=\R\cup\{\infty\}$ 
({\it i.e.} with the compactification at $+\infty$) leaving $\infty$ fixed and let ${\bar \R}\times\R$ 
be the action groupoid. The space of leaves ${\cal L}$ can be characterized as the subgroupoid 
$\{(s,t) | ~~ s\geq 0,~ s+t\geq 0\}\cup\{(\infty,0)\}$ over $\overline{\R_{\geq0}}=\R_{\geq 0}\cup \{\infty\}$. 
In fact the lagrangian leaves are the contour levels of the groupoid morphism $(f,F)$, where 
$f:\sphere^2\rightarrow \overline{\R_{\geq 0}}$ and $F:\G(\sphere^2,\pi)\rightarrow {\cal L}$ can be 
read in the symplectic chart as:
\begin{equation}\label{gpd_hom}
f(z)=\log(1+|z|^2) ~~,~~~~ F(x,y)=(\log (1+|x|^2),\log(\frac{1+|y|^2}{1+|x|^2}))  \;.
\end{equation}

\smallskip
\begin{rem}{\rm ({\it The polarization is singular}).
The point $\infty\in{\cal L}$ is a singular point from the point of view of the smooth structure, 
so that ${\cal L}$ inherits from the ambient space just the topology. The tangent spaces to the contour 
levels of $F$ fail to define a continuous distribution of ${\cal G}(\sphere^2,\pi)$, even in the generalized sense. 
Let us check it by computing the coordinate vector fields for $w\not=0$ as
\begin{eqnarray*}
\frac\partial{\partial y}&=& -\frac{w}{\bar{w}}\frac{1}{1+|w|^2}\frac\partial{\partial \bar p_S}\,,\\
\frac\partial{\partial \bar y}&=& -\frac{\bar{w}}{w}\frac{1}{1+|w|^2}\frac\partial{\partial p_S}\,,\\
\frac\partial{\partial x}&=& -w^2\frac\partial{\partial w}+\frac{w}{\bar{w}}\ \frac{1-\bar{w}\bar{p}_S}{1+|w|^2}
\frac\partial{\partial \bar p_S}+
\frac{w\, p_S}{1+|w|^2}(1+2|w|^2)\frac\partial{\partial p_S}\,,\\
\frac\partial{\partial \bar{x}}&=& -\bar{w}^2\frac\partial{\partial \bar{w}}+\frac{\bar{w}}{w}\ \frac{1-wp_S}
{1+|w|^2}\frac\partial{\partial p_S}+
\frac{\bar{w}\, \bar{p}_S}{1+|w|^2}(1+2|w|^2)\frac\partial{\partial \bar{p}_S}\,.
\end{eqnarray*}

The polarization in $w\not=0$ is spanned by the Hamiltonian fields
$$
\chi_+=i(1+|x|^2)(\bar x\partial_{\bar x}-x\partial_x),\qquad
\chi_-=i(1+|y|^2)(\bar y\partial_{\bar y}-y\partial_y)\,,
$$
that in terms of the coordinates $w$ and $p_S$ read
\begin{eqnarray*}
\chi_+ = i(1+|w|^2)(-\frac{1}{w}\frac{\partial}{\partial \bar{w}}+ \frac{1}{\bar{w}}\frac{\partial}{\partial w})&+& 
i \frac{1}{|w|^2\bar{w}}\left[-1+2\bar{w}\bar{p}_S(1+|w|^2)\right]\frac{\partial}{\partial \bar{p}_S}\cr
&-&i \frac{1}{|w|^2w}\left[-1+ 2wp_S(1+|w|^2)\right]\frac{\partial}{\partial p_S}
\end{eqnarray*}
$$
\chi_- = i A(w,p_S) \left[ \frac{1}{w\bar{w}^2}(1+w\delta y)\frac{\partial}{\partial\bar{p}_S} - 
\frac{1}{\bar{w}w^2}(1+\bar{w}\delta \bar{y})\frac{\partial}{\partial p_S}\right]
$$
where
$$
A(w,p) = (1+|y|^2) \frac{|w|^2}{1+|w|^2} \;,\;\;\;\;\;    \delta y = y-\frac{1}{w}\ .
$$

The polarization in $(w,0)$ is spanned by
$$
\chi_++\chi_-=i\ \frac{(1+|w|^2)}{|w|^2}\left(w\frac{\partial}
{\partial w}-\bar{w}\frac{\partial}{\partial \bar{w}}\right)\,,~~~~~ \chi_-=
\frac{i}{|w|^2}\left(\frac{1}{\bar{w}}\frac{\partial}{\partial\bar{p}}-\frac{1}{w}
\frac{\partial}{\partial p}\right)\;;
$$
it is clear that no regular vector fields survive in the limit $w\rightarrow 0$}. $\square$ \end{rem}

Let us analyze {\it Bohr-Sommerfeld} conditions of the Lagrangian leaves. A Lagrangian leaf is a Bohr-Sommerfeld 
leaf if the flat connection obtained by restriction of the prequantization has trivial  holonomy. Bohr-Sommerfeld 
condition means that both $\int_{|x|^2=F_+} \Theta$ and $\int_{|y|^2=F_-} \Theta$ live in $2\pi\hbar \Z$ where 
$\Theta$ is any primitive of the symplectic form (\ref{symplectic_form_grpd_pdls_sphere}). By direct computation 
this means that the leaf defined by $(F_+,F_-)$ is Bohr-Sommerfeld if and only if there exist $n_\pm\in\Z_{\geq 0}$ 
such that
$$
F_{\pm} = e^{\hbar n_\pm} - 1 \;.
$$
Let us describe the $l$ (or $r$) image of Bohr-Sommerfeld leaves on $\sphere^2$. Let $\alpha=z/(1+|z|^2)$, 
$\tau=1/(1+|z|^2)$ be the usual parametrization of $\sphere^2$. The Bohr-Sommerfeld leaves on $\sphere^2$ are 
given by
$$
\tau= e^{-\hbar n} \;, ~~~~~ n\geq 0\;,
$$
{\it i.e.} the spectrum of $\rho(\tau)$ given in (\ref{unitary_rep_podles_sphere}) for $q=e^{-\hbar/2}$. 
This formula differs by a factor from (\ref{ellestartau}) obtained by means of a complex polarization.
Indeed Bohr-Sommerfeld rules do not take into account the metaplectic correction as in the energy 
representation of the one dimensional harmonic oscillator.

\begin{rem}{\rm
According to (\ref{gpd_hom}), the Bohr-Sommerfeld leaves identify the subgroupoid of ${\cal L}$ 
described as ${\cal L}_{BS}=\{(\hbar m,\hbar n ),~ m+n\geq 0, m\geq 0\}\cup \{\infty\}$ with space of 
units given by $(\overline{\R_{\geq 0}})_{BS}=\{ \hbar n, n\geq 0\}\cup\{\infty\}$. Equipped with the relative 
topology, this groupoid of BS--leaves coincide with Sheu's groupoid $\G_S$ described in Section 
\ref{podles_sphere}. }
\end{rem}

\begin{rem}{\rm
We see that the modular function $c_{V_{\sphere^2}}$, computed in (\ref{symplectic_modular_function}), descends 
to a cocycle of $\G_S$ that reads $c_{V_{\sphere^2}}=\hbar c_1 $, where $c_1$ is the cocycle introduced in Section 
\ref{ToplitzGrpd}. The corresponding modular automorphism $A_{c_1}(-i \hbar)$ is exactly the algebraic van den 
Bergh automorphism (which, for the standard Podle\`s sphere was computed in \cite{Kra}).
This relation is not surprising since, as shown by Dolgushev in \cite{Dol}, in the context of formal deformation 
quantization the modular vector field is quantized by a derivation of the deformed algebra which exponentiate to 
the van den Bergh automorphism. Geometric quantization through symplectic groupoids confirms, through a different 
path, the same result (let us recall that Dolgushev's theorems do no apply directly here). The Poisson modular
class can thus be shown to give rise to two fundamental objects arising in quantization: the modular operator of 
Tomita--Takesaki theory on the analytic side and the van den Bergh dualizing bimodule on the algebraic side.}
\end{rem}

\bigskip

\subsection{A complex polarization}\label{complex_polarization}

We introduce here a complex polarization. Since the symplectic groupoid restricted to the symplectic leaf is a 
K\"ahler manifold, we can take the K\"ahler polarization and then extend it with the vertical polarization on the 
zero dimensional leaf. We will show that we get a smooth lagrangian distribution.

Let the polarization $F$ be defined as $F=\langle\frac{\partial}{\partial x},\frac{\partial}{\partial \bar{y}}
\rangle$ for finite $(x,y)$ and $F=(T^*_{N}\sphere^2)_\C$ in $(w=0,p_S)$. One can check that this is the
polarization considered in Theorem 8.1 of \cite{Hawkins} for a general K\"ahler-Poisson manifold.

\begin{lemma}
\label{smoothness_polarization}
$F$ is a positive and multiplicative polarization. Moreover $\det F$ is a trivial line bundle.
\end{lemma}
{\it Proof}. For any $(w,p)$ let us compute
\begin{eqnarray*}
v_1(w,p_S) &\equiv& \frac{\bar{w}}{w} \frac{\partial}{\partial x} = -|w|^2\frac{\partial}{\partial w} + 
\frac{1-\bar{w}\bar{p}_S}{1+|w|^2}\frac{\partial}{\partial\bar{p}_S} + \frac{\bar{w}p_S}{1+|w|^2}(1+2|w|^2)
\frac{\partial}{\partial p_S} \cr v_2(w,p_S) &\equiv & \frac{w}{\bar{w}}\frac{\partial}{\partial \bar{y}}= 
-\frac{1}{1+|w|^2}\frac{\partial}{\partial p_S}\;.
\end{eqnarray*}
Since $v_1(0,p_S)=\partial/\partial \bar{p}_S$ and $v_2(0,p_S)=\partial/\partial p_S$, $\{v_1,v_2\}$ is a 
local smooth basis for $F$ and $\partial_x\wedge\partial_{\bar{y}}=v_1\wedge v_2$ is a global non vanishing 
section for $\det F$. Multiplicativity and positivity are direct checks. $\square$

\medskip
\begin{rem}{\rm
The associated real distributions are singular; in fact ${\cal D}=F\cap \bar{F}$ is $0$ on finite $(x,y)$ 
and $(T^*_N\sphere^2)_\C$ in $(w=0,p_S)$}.
\end{rem}
\smallskip
A basis of hamiltonian vector fields is $\{\chi_{\bar x},\chi_{y}\}$; by construction they are covariantly constant 
under the Bott connection. We identify $F$ with $F^\perp=\Omega(F)$, so that the basis of covariantly constant 
sections is $\underline{b}=\{d\bar{x},dy\}$. Since $\det F^\perp$ is trivial, its unique square root 
$\sqrt{\det F^\perp}$ is the trivial line bundle.

Let us choose the following local primitive of the symplectic form (\ref{symplectic_form_grpd_pdls_sphere})
\begin{equation}
\label{primitive}
\Theta= \frac{1}{2i}\log(1+|x|^2)\left(\frac{d{\bar x}}{{\bar x}}-\frac{dx}{x}\right) -\frac{1}{2i}\log(1+|y|^2)
\left(\frac{d{\bar y}}{{\bar y}}-\frac{dy}{y}\right)\;.
\end{equation}
The local primitive $\Theta$ is multiplicative, {\it i.e.} $\partial^*\Theta=0$ so that the corresponding 
prequantization cocycle is $1$. We stress the fact that this is true only locally, in fact the obstruction 
to have a global multiplicative primitive has been discussed at the beginning of this Section.

A polarized section $\sigma\otimes\sqrt{\underline{b}}\in C^\infty(T^*\sphere^2)\otimes 
C^\infty(\sqrt{\det F^\perp})$ satisfies
$$
\partial_x \sigma+\frac{1}{2\hbar}\frac{\log(1+|x|^2)}{x}\sigma= \partial_{\bar y} \sigma+\frac{1}{2\hbar}
\frac{\log(1+|y|^2)}{\bar{y}}\sigma = 0\;.
$$
Observe that, since $\underline{b}$ and $\Theta$ are only local, $\sigma$ is not a global function on 
$T^*\sphere^2$. If one uses the global section given in the proof of Lemma (\ref{smoothness_polarization}), 
then he has to add the contribution of the Bott connection. The solution can be written as
\begin{equation}
\label{polarized_sections}
\sigma= \psi(\bar{x},y) e^{\frac{1}{2\hbar}(Li_2(-|x|^2)+Li_2(-|y|^2)})\;,~~~
\end{equation}
where $ Li_2(t)=-\int_0^{t} d\tau \frac{\log(1-\tau)}{\tau}$ is the integral expression for the dilogarithm. 
According to (\ref{symplectic_isomorphism}), the scalar product reads
\begin{equation}
\label{scalar_product}
\langle\sigma_1,\sigma_2 \rangle = \int_{\C^2} d^2x d^2y \ e^{\frac{1}{\hbar}(Li_2(-|x|^2)+Li_2(-|y|^2))} 
\frac{\overline{\psi_1(\bar{x},y)} \psi_2(\bar{x},y)}{\sqrt{1+|x|^2}\sqrt{1+|y|^2}}\;,
\end{equation}
where $d^2x=dx d\bar{x}/i$.
The (local) observable $f=\log(1+|x|^2)$ is quantizable; indeed its hamiltonian vector field $\chi_f=i(\bar{x}
\partial_{\bar x}-x\partial_x)$ satisfies 
$[\chi_f,\chi_{\bar x}]=i\chi_{\bar x}$. So the quantization rule reads (see (\ref{quantization}))
$$
\hat{f}\psi = \hbar \bar{x}\frac{\partial}{\partial \bar{x}}\psi+\frac{\hbar}{2}\psi\;.
$$

\begin{lemma}
The states $\sigma_{m,n}=\bar{x}^my^n\exp(Li_2(-|x|^2)+Li_2(-|y|^2))/2\hbar$ are orthogonal and normalizable.
\end{lemma}
{\it Proof}. The norm of $\sigma_{m,n}$ is computed by using the scalar product (\ref{scalar_product}) as 
$||\sigma_{m,n}||^2=A_mA_n$, where
$$ A_m= 2\pi \int_0^\infty dt\ \frac{t^m}{\sqrt{1+t}} e^{\frac{1}{\hbar} Li_2(-t)} ~~.$$ 
We need the following functional relation valid for $t>0$ (see \cite{maximon})
\begin{equation}\label{dilogarithm}
Li_2(-t) = - Li_2(-\frac{1}{t}) - \frac{\pi^2}{6} -\frac{1}{2} (\log t)^2 \;,
\end{equation}
that implies the asymptotic behaviour for $t\gg 1$. As a consequence, for $t>t_0$, $Li_2(-t)\leq C_{t_0} - 
1/2 (\log t)^2$. Since $\exp(-(\log t)^2/2\hbar)$ is rapidly decreasing then $A_m$ is bounded by a converging 
integral. 
$\square$
\smallskip

The eigenvalue problem for $\hat{f}$ is the same as the usual two dimensional harmonic oscillator (with half zero 
point energy), so that, for any $m\in\N$, 
$$
\hat{f}\sigma_{m,n} = \hbar (m+\frac{1}{2})\sigma_{m,n}\;.
$$
Remark that $\tau=1/(1+|x|^2)$ is a global function on $\sphere^2$ that can be embedded as $l^*\tau=e^{-f}$ in 
$C^\infty(T^*\sphere^2)$, so that 
\begin{equation}
\label{ellestartau}
\widehat{l^*\tau}\sigma_{m,n} = e^{-\hbar(m+1/2)} \sigma_{m,n}\;.
\end{equation}
Moreover the modular function (\ref{symplectic_modular_function}) can be written as 
$D_{V_{\sphere^2}}=l^*\tau/r^*\tau$ so that it is quantizable and we get
$$
\Delta \sigma_{m,n}=\widehat{D}_{V_{\sphere^2}}\sigma_{m,n} = e^{-\hbar(m-n)} \sigma_{m,n}
$$
that we recognize as the modular operator of $C^*({\cal O}_1)$ described in Subsection \ref{ToplitzGrpd}.

\bigskip
\bigskip

\section{The left Hilbert algebra from polarized sections}
We want now to define a convolution product between the polarized sections of the complex polarization studied in 
Subsection \ref{complex_polarization}. Our task is to define a left Hilbert algebra by using the scalar product 
(\ref{scalar_product}) and compute the corresponding modular automorphism.

If $\sigma_i$, $i=1,2$, are polarized sections,  then $\sigma_1(x,z)\sigma_2(z,y)$ lives in the square root of 
$\det F^\perp_{(x,z)}\otimes\det F^\perp_{(z,y)}=$ $\det F^\perp_{(x,y)}\otimes\det T^*_{(x,z)}l^{-1}(x)$. Since 
we want to integrate on $l^{-1}(x)$, we need to tensor it with a nonvanishing section of (the square root of) 
$\det T^*_{(x,z)}l^{-1}(x)$ in order to get a top form. This can be done with $\Omega|_{l^{-1}(x)}$, so that we 
get the following definition for the convolution product and involution of polarized sections $\sigma_i$, 
\begin{equation}
\label{convolution_from_geometric_q}
\sigma_1 *\sigma_2(x,y) = \int_{\C} d^2z \frac{1}{\sqrt{1+|z|^2}} \sigma_1(x,z)\sigma_2(z,y)\;, ~~~
\sigma^*(x,y) = \overline{\sigma(y,x)}\;. 
\end{equation}

\begin{rem}
{\rm Formula (\ref{convolution_from_geometric_q}) coincides with the prescription given in Section 5.3-5.4 of 
\cite{Hawkins}. Indeed we have that $\Omega_F=\det F^\perp$, where $\Omega_F$ is defined in (5.4) of \cite{Hawkins}, 
and the above steps leading to the convolution formula are the same as applying Lemma 5.3 of \cite{Hawkins}.}
\end{rem}

Equipped with its scalar product (\ref{scalar_product}), the space of polarized sections becomes a left Hilbert
algebra, as one can check directly. The modular operator is defined from the polar decomposition of the involution
operator $S(\sigma)=\sigma^*$, {\it i.e.} $\Delta=S^\dagger S$. One immediately checks that $S^\dagger=S$ so that $\Delta=1$ !
Indeed $S^\dagger \sigma_{m,n} = S \sigma_{m,n} = \sigma_{n,m}$, since the symplectic scalar product
(\ref{scalar_product}) is symmetric in the exchange of $x$ with $y$. 

Both scalar product and convolution formula depend on choices that have been done by using uniquely the symplectic
structures. We are going to show that the right choice is done instead by using the left Haar system introduced in 
(\ref{haar_system_podles}).  In fact, the Haar system (\ref{haar_system_podles}) consists in giving  a volume form 
on $l^{-1}(x)$ for any $x\in\sphere^2$, depending on a non vanishing $\Lambda\in C^\infty(\sphere^2)$, that can be used to 
define the convolution algebra; with this choice we get
\begin{equation}
\label{convolution_from_groupoid}
\sigma_1 *_\Lambda\sigma_2(x,y) = \int_{\C} d^2z \sqrt{\Lambda(z)} \sigma_1(x,z)\sigma_2(z,y)\;, ~~~
\sigma^*(x,y) = \overline{\sigma(y,x)}\;. 
\end{equation}
  
The scalar product (\ref{scalar_product}) depends on the choice of the isomorphism (\ref{symplectic_isomorphism}), 
that in our case is a trivialization of $\det T^*(T^*\sphere^2)$, that in the purely symplectic case is naturally 
done with the symplectic volume. If one fixes a Haar system on the groupoid, then other choices are possible. 
In fact, for any choice of the volume form $V_{\sphere^2,\rho}=\rho V_{\sphere^2}$ on $\sphere^2$, for 
$V_{\sphere^2}$ being the round volume form and some positive $\rho\in C^\infty(\sphere^2)$, there is defined the 
following volume form on $T^*\sphere^2$
$$\nu_{(\rho,\Lambda)}= \rho(x)\Lambda(y)D_{V_{\sphere^2}}\Omega^2\;.$$ 
Remark that this procedure defines a subset of all possible volume forms and that the symplectic volume $\Omega^2$ 
is not in this set, due to the non triviality of the modular class. If we use $\nu_{(\rho,\Lambda)}^{-1}$ to 
trivialize $\det T^*(T^*\sphere^2)$ we get for the scalar product 
\begin{equation}\label{groupoid_scalar_product}
\langle \sigma_1,\sigma_2\rangle_{(\rho,\Lambda)} = \int_{\C} d^2x d^2y \frac{\sqrt{\rho(y)\Lambda(x)
D_{V_{\sphere^2}}^{-1}(x,y)}}{\sqrt{(1+|x|^2)(1+|y|^2)}}\  \overline{\sigma_1(x,y)}\sigma_2(x,y)\;;
\end{equation}
moreover, if we choose $\Lambda(y)=\Lambda(|y|^2)$ and $\rho(x)=\rho(|x|^2)$, then 
$||\sigma_{m,n}||^2_{(\lambda,\Lambda)}= \ell_m r_n$, where
\begin{equation}\label{integrals}
\ell_m = 2\pi\int_0^\infty dt\ t^m \sqrt{\Lambda(t)} e^{\frac{1}{\hbar}Li_2(-t)}\,,~~~~~
r_m = 2\pi \int_0^\infty dt \frac{t^m}{1+t} \sqrt{\rho(t)} e^{\frac{1}{\hbar}Li_2(-t)}
\,.
\end{equation}
It is important to state the following asymptotic behaviour.

\begin{lemma}
\label{asymptotic}
The asymptotic expansion of $\ell_n/r_n$ for large $n$ and fixed $\hbar$ is
$$
\frac{r_n}{\ell_n} \sim e^{-\hbar(n+\frac{1}{2})} \sqrt{\frac{\rho(\infty)}{\Lambda(\infty)}} \;.
$$ 
\end{lemma}
{\it Proof}. Both $\ell_n$ and $r_n$ are divergent in this asymptotic limit, as we will check later, so that
we can subtract from the integrals in (\ref{integrals}) the contribution on $(0,1)$ which are finite. By 
using (\ref{dilogarithm}), we see that we have to evaluate the asymptotic of integrals of this form
$$
I_n(F) = 2\pi \int_1^\infty dt \  e^{n(\log t -\frac{1}{2\hbar}(\log t)^2)} F(t) \;,
$$
where $r_n = I_{n-1}(F_r)$, $\ell_n= I_n(F_\ell)$, $F_\ell(t)=\exp[(-Li_2(-1/t)-\pi^2/6)/\hbar]\ 
\sqrt{\Lambda(t)}$ and $F_r(t)=\exp[(-Li_2(-1/t)-\pi^2/6)/\hbar]\ \sqrt{\rho(t)}/(1+1/t)$. 

With $s=\log t/(\hbar(n+1))$ we get
$$
I_n(F) = 2\pi \hbar (n+1) e^{\hbar(n+1)^2/2}\int_0^\infty ds\ e^{-\frac{\hbar}{2}(n+1)^2(s-1)^2} F(e^{\hbar(n+1)s})\;. 
$$
Since both $F_\ell$ and $F_r$ are bounded in the integration domain, we get the asymptotic of $I_n(F)$ by the 
saddle point evaluation. We collect so
$$
r_n \sim 2\pi \sqrt{2\pi\hbar}\ e^{-\pi^2/6\hbar}\ e^{\hbar n^2/2}~,~~~~ \ell_n \sim 2\pi \sqrt{2\pi\hbar}\ 
e^{-\pi^2/6\hbar}\ e^{\hbar(n+1)^2/2}~, 
$$
from which the result follows. $\square$
\smallskip

Recalling now the discussion in Subsection \ref{ToplitzGrpd} about the GNS construction associated to the quasi 
invariant measure $\mu_\hbar$ with modular function $e^{\hbar c_1}$, where $c_1(m,n)=n$, we state the following 
Proposition.

\smallskip
\begin{prop}
\label{left_hilbert_algebra}
The vector space ${\cal A}_{\widehat{D}}$ spanned by the eigenvectors of $\widehat{D}_{V_{\sphere^2}}=\Delta$ equipped 
with the convolution and involution {\rm (\ref{convolution_from_groupoid})} and the scalar product 
{\rm (\ref{groupoid_scalar_product})} is a left Hilbert algebra. The map $e:{\cal A}_{\widehat D}\rightarrow 
C_c({\cal O}_1)$ defined by
\begin{equation}\label{hilb_alg_isom}
 e(\sigma_{m,n})= e_{m,n-m}\sqrt{\ell_m\ell_{n}} 
\end{equation}
is an algebra homomorphism and an isometry if $C_c(\O_1)$ is equipped with the scalar product of the GNS state 
$\phi_{\mu_{(\rho,\Lambda)}}$ defined by the quasi invariant probability measure
$$
\mu_{(\rho,\Lambda)}(m) = \frac{r_m}{\ell_m}\;.
$$ 
This measure is equivalent to $\mu_\hbar$ so that the modular automorphism is cohomologous to $c_1$, {\it i.e.} 
$c_{(\rho,\Lambda)} = c_1 + \partial^*\varphi $, where $\varphi\in C_c(\overline{\N})$ is defined as
$\varphi(m)=\frac{1}{\hbar}\log (r_m/\ell_m) +m$ and $\varphi(\infty)=-\frac{1}{2}+\frac{1}{2}\log
\frac{\rho(\infty)}{\Lambda(\infty)}$. 
\end{prop}
{\it Proof.} By direct computation one shows that (\ref{hilb_alg_isom}) respects the relations 
(\ref{cuntz_grpd_relations}) of $C_c(\O_1)$ and  $$\phi_{\mu_{(\rho,\Lambda)}}(e(\sigma_{m,n})^**
e(\sigma_{m,n}))=\phi_{\mu_{(\rho,\Lambda)}}(e_{n,0})\ell_m\ell_n = \mu_{(\rho,\Lambda)}(n)\ell_m \ell_n = 
||\sigma_{m,n}||^2_{(\rho,\Lambda)}.$$ 

Since ${\cal A}_{\widehat D}^2={\cal A}_{\widehat D}$ and $e({\cal A}_{\widehat D})$ is dense in 
${\cal H}_{\mu(\rho,\Lambda)}$ it follows that ${\cal A}_{\widehat D}$ is an Hilbert algebra.
The eigenvectors of the modular operator  $\Delta_{(\rho,\Lambda)}=S^\dagger S$ are $\sigma_{m,n}$ with 
eigenvalue $\ell_n r_m/\ell_m r_n$ and the formula for the modular cocycle $c_{(\rho,\Lambda)}=
\log D_{(\rho,\Lambda)}/\hbar$ follows. Moreover, $\lim_{m\rightarrow\infty}\varphi(m)=\varphi(\infty)$ 
as a consequence of Lemma \ref{asymptotic}. $\square$ 

\medskip

\begin{rem}{\rm \begin{itemize} 
\item[$i$)] The choice of $(\rho,\Lambda)$ is not relevant, since it doesn't change the class of the 
modular cocycle. On the contrary, no choice of $(\rho,\Lambda)$ in (\ref{groupoid_scalar_product}) can 
recover the symplectic scalar product (\ref{scalar_product}).
\item[$ii$)] The asymptotic behaviour stated in Lemma \ref{asymptotic} is the condition that assures that 
the identity $1_\Lambda=\sum_{n\geq 0} e_{n,0}$ has finite norm with respect to the scalar product 
(\ref{groupoid_scalar_product}). 
\end{itemize}
}
\end{rem}

\bigskip
\bigskip

\section{Conclusions}

\noindent {\bf Von Neumann versus $C^*$-algebra}. 
The outcome of our construction is an Hilbert algebra and, by applying a standard procedure, 
a Von Neumann algebra. This framework is very natural from the quantization point of view and in the paradigm of 
non commutative geometry it means that we are looking to the quantum space from the point of view of measure theory. 
What one should do in order to get the right topology, {\it i.e.} the right $C^*$-algebra ? The closure of 
${\cal A}_{\widehat D}$ with respect to the operator norm gives the $C^*$-algebra ${\cal K}$ of compact operators. 
We know that the correct answer is obtained by adjoining the identity.

\smallskip 
\noindent {\bf Role of characteristic classes}. In the cohomology of a Poisson manifold there are two distinct 
classes that play an important role in the quantization procedure: the class of the Poisson tensor itself and 
the modular class. 

In the example we analyzed in this paper, they are both non trivial. The non vanishing of the class of the Poisson 
tensor is responsible for the appearing of a non trivial prequantization two-cocycle that twists the convolution 
algebra of sections. In this example the class is trivial once we remove a measure zero set from the symplectic 
groupoid, and for this reason we didn't need the explicit formula of the prequantization cocycle. In more general 
situations we expect that an explicit formula is unavoidable and moreover that a compatibility between the 
polarization and the cocycle must be imposed.

The modular class plays a major role in our construction. First it puts restrictions on the choice of 
polarization. With such choices the modular vector field is quantized to the modular operator which describes an 
intrinsic property of the quantum measurable space. In fact it plays a relevant role in defining the scalar product 
between polarized sections. This construction can be seen as the analytic analogue of the algebraic quantization 
of the modular class provided by Dolgushev (see \cite{Dol}).

\smallskip
\noindent {\bf Further examples}. If one is interested in understanding the procedure of quantization we believe 
that our analysis shows the importance of studying concrete examples. It is very natural to try to extend the 
results of this paper to compact quantum groups, mainly $SU_q(n)$, and quantum homogeneous spaces, for example 
$\C P_q(n)$ or other quantum spheres. Sheu's papers showed that their $C^*$-algebras can be described as groupoid 
$C^*$-algebras so that we can expect that we can obtain results similar to those described here. Some difficulties 
that here we avoided due to the simplicity of the example will appear: for instance an explicit expression of the 
quantization cocycle will be probably needed. It has to be expected, though, that real singular polarizations will 
be needed in these cases as well; a general theory of such polarizations is at present still missing.

One can also consider other Poisson structures with the same foliation than the standard Podle\`s sphere but with a 
different degree of singularity, for instance the quartic one. As Poisson manifolds, they are completely inequivalent 
(for instance infinite dimensional Poisson cohomology versus finite dimensional). The one with quartic singularity 
carries an abelian action and can be therefore explicitly integrated and quantized in this scheme. The topology of the 
quantum space should be the same as the Podle\`s sphere we are considering. It would be relevant to understand how 
to describe additional structures (like the smooth structure) in order to distinguish the quantum spaces.

\bigskip
\bigskip

\noindent{\bf\Large Acknowledgement}:
\bigskip

\noindent We are grateful to Domenico Seminara and Filippo Colomo for useful discussions.

\eject

\end{document}